\documentclass{amsart}

\usepackage{qayum_thesis}

\title[Calculation of $\UNil$ for $C_2$]{Calculation of UNil for the\\ cyclic group of order two}
\author{Qayum Khan}
\address{Department of Mathematics\\ Vanderbilt University\\ Nashville, TN 37240 U.S.A.}
\email{qayum.khan@vanderbilt.edu}

\ncm{\map}{map} \ncm{\mor}{mor} \ncm{\Or}{Or} \ncc{\cF}{F} \ncq{\qF}{F} \ncc{\cG}{G}
\ncm{\cyc}{cyc} \ncm{\dih}{dih} \ncm{\fbc}{fbc} \ncm{\fin}{fin} \ncm{\mcyc}{mcyc}
\ncm{\mdih}{mdih} \ncm{\pfin}{pfin} \ncm{\vc}{vc}

\begin{document}

\begin{abstract}
Cappell's unitary nilpotent groups $\UNil_*^h(R;R,R)$ are calculated
for the integral group ring $R = \mathbb{Z}[C_2]$ of the cyclic
group $C_2$ of order two. Specifically, they are determined as
modules over the Verschiebung algebra $\cV$ using the
Connolly--Ranicki isomorphism \cite{CR} and the Connolly--Davis
relations \cite{CD}.
\end{abstract}

\maketitle

\section{Introduction}

Consider the simplest nontrivial 2-group
\[
C_2 = \gens{T\ST T^2=1}. \] Observe that the integral group ring
$\Z[C_2]$ fits into Rim's cartesian square
\[\begin{CD}
\Z[C_2] @>{i^-}>> & \Z\\
@V{i^+}VV & @V{j^-}VV\\
\Z @>{j^+}>> & \F_2
\end{CD}\]
of rings with involution, where $i^\pm(T) = \pm 1$.  We focus on
piecing together its $\UNil$ as a module from the $\UNil$ of the
component rings $\Z$ and $\F_2$, using a Mayer--Vietoris sequence in
algebraic $L$-theory. The additional structure\footnote{The
$\cV$-module structure is induced by the Connolly--Ranicki
isomorphism (Thm. \ref{Thm_CRiso}).} on
\[
\UNil_*^h(R) := \UNil_*^h(R;R,R) \iso NL_*^h(R)
\]
computed below is its covariant (pushforward) module structure over
the \textbf{Verschiebung algebra}
\[
\cV := \Z[V_n \ST n>0 ] = \Z[{V_p \ST p \text{ prime}}]
 \]
of $n$-th power operators
\[
V_n := (x \mapsto x^n)
\]
on polynomial rings $R[x]$. An analogous structure in algebraic
$K$-theory has been studied by Joachim Grunewald \cite[\S
4.2.2]{GrunewaldLocalization} for the Bass $\Nil$-groups
\[
\wt{\Nil}_*(R) = NK_{*+1}(R).
\]

The following two theorems are the main results of this paper. In
the instance $F=C_2$, since the lower $\wt{\Nil}_i$- and
$NK_{i+1}$-groups vanish \cite{Harmon}, we may replace the
$\gens{-\infty}$ decoration with the $h$ decoration. The first main
theorem provides a general vanishing result and a classifying
isomorphism, specializing \cite{Khan_Sylow2NormalAbelian}.

\begin{thm}\label{Thm_VanishExponentTwo}
Suppose $F$ is a finite group that contains a normal Sylow
2-subgroup of exponent two.  If $n \equiv 0,1 \pmod{4}$, then the
following abelian group vanishes:
\[
\UNil_n^{\gens{-\infty}}(\Z[F]) = 0.
\]
Furthermore, if $n \equiv 2 \pmod{4}$, then the following induced
map is an isomorphism:
\[
\UNil_n^{\gens{-\infty}}(\Z[F]) \xra{\quad}
\UNil_n^{\gens{-\infty}}(\F_2) \xra{r\iso}
NL_n^{\gens{-\infty}}(\F_2) \xra{\Arf \iso} x\F_2[x]/(f^2-f).
\]
\end{thm}

The second main theorem examines non-vanishing in the remaining
dimensions.

\begin{thm}\label{Thm_OrderTwoAnswer}
If $n \equiv 3 \pmod{4}$, then there exists a decomposition
\[
\UNil_n^{\gens{-\infty}}(\Z[C_2]) \iso
\UNil_{n+1}^{\gens{-\infty}}(\F_2) \oplus
\UNil_n^{\gens{-\infty}}(\Z) \oplus \UNil_n^{\gens{-\infty}}(\Z).
\]
\end{thm}

\begin{proof}
Immediate from Theorems \ref{Thm_OrderTwo} and \ref{Thm_NL3}.
\end{proof}

\section{Definitions, relations, and decompositions}

Unless specified otherwise, all the surgery groups $L, \UNil, NL$ in
in this paper shall have the $h$ decoration with respect to the
algebraic $K$-groups $\wt{K}_1, \wt{\Nil}_0, NK_1$.

\begin{defn}[Bass]
Let $R$ be a ring with involution.  For each $n \in \Z$, define the
abelian group
\[
NL_n(R) := \Ker\prn{\aug_0: L_n(R[x]) \longra L_n(R)}.
\]
Therefore there is a natural decomposition
\[
L_n(R[x]) = L_n(R) \oplus NL_n(R).
\]
\end{defn}

The subsequent statements are technical tools for the main theorems.
The first is the Connolly--Ranicki isomorphism \cite[Thm. A]{CR},
which is a fundamental equivalence in the computation of a certain
class of $\UNil$-groups.

\begin{thm}[Connolly--Ranicki]\label{Thm_CRiso}
Let $R$ be a ring with involution. Then for all $n \in \Z$, there is
a natural isomorphism
\[
r^h: \UNil_n^h(R;R,R) \longra NL_n^h(R)
\]
which descends to a natural isomorphism
\[
r^{\gens{-\infty}}: \UNil_n^{\gens{-\infty}}(R;R,R) \longra
NL_n^{\gens{-\infty}}(R).
\]
\end{thm}

\begin{rem}\label{Rem_NLcomp}
According to Connolly--Ko\'zniewski \cite{CK}, Connolly--Ranicki
\cite{CR}, and Connolly--Davis \cite{CD}, the group
$NL_{\odd}(\F_2)$ vanishes, and the Arf invariant is an isomorphism:
\[
\Arf: NL_{\even}(\F_2) \xra{\quad} \frac{x \F_2[x]}{(f^2-f)}.
\]
The inverse of $\Arf$ is given by the map $q \longmapsto P_{q,1}$,
where for all $p,g \in \Z[x]$ the symplectic form $P_{p,g}$ is
defined by
\[
P_{p,g} := \prn{\bigoplus_2 \F_2[x], \PMatrix{0 & 1\\ 1 & 0},
\PMatrix{p\\ g}}.
\]
Also the group $NL_n(\Z)$ vanishes if $n \equiv 0,1 \pmod{4}$, the
induced map to $NL_2(\F_2)$ is an isomorphism if $n \equiv 2
\pmod{4}$, and there is a two-stage obstruction theory \cite[Proof
1.7]{CD} if $n \equiv 3 \pmod{4}$:
\[
0 \xra{\quad} \frac{x \F_2[x]}{(f^2-f)} \xra{\en \cP \en} NL_3(\Z)
\xra{\en B\en} x \F_2[x] \x x \F_2[x] \xra{\quad} 0.
\]
It is given primarily by certain characteristic numbers $B$ in Wu
classes of $(-1)$-quadratic linking forms over $(\Z[x],2)$, and
secondarily by the Arf invariant, of even linking forms $\cP$, over
the function field $\F_2(x)$.
\end{rem}

\begin{thm}\label{Thm_OrderTwo}
Consider $P = C_2$ with trivial orientation character. Then, as
Verschiebung modules, there is a decomposition
\[
NL_3(\Z[C_2]) = NL_3(\Z) \oplus \wt{NL}_3(\Z[C_2])
\]
and there is an exact sequence (constituting a three-stage
obstruction theory):
\[\begin{CD}
0 @>>> & NL_0(\F_2) @>{\wt{\bdry}}>> & \wt{NL}_3(\Z[C_2]) @>{i^-}>>
NL_3(\Z) @>>> 0.
\end{CD}\]
\end{thm}

Ingredients for the next theorem are as follows. Given a ring $A$
with involution and $\epsilon = \pm 1$, there is an identification
\cite[Prop. 1.6.4]{RanickiExact} between split $\epsilon$-quadratic
formations over $A$ and connected 1-dimensional $\epsilon$-quadratic
complexes over $A$.  The identification between
$(-\epsilon)$-quadratic linking forms over $(A,(2)^\infty)$ and
resolutions by $(2)^\infty$-acyclic 1-dimensional
$\epsilon$-quadratic complexes over $A$ is given by
\cite[Proposition 3.4.1]{RanickiExact}.

The determination of the above extension (\ref{Thm_OrderTwo}) of
abelian groups involves algebraic gluing of quadratic complexes
\cite[\S 1.7]{RanickiExact}, given below (\ref{Def_QM}) by a choice
$\cM$ of set-wise section.  Recall from group cohomology that an
extension of abelian groups
\[
0 \xra{\quad} A \xra{\quad} B \xra{\quad} C \xra{\quad} 0
\]
and a choice of set-wise section $s: C \to B$ determine a
\textbf{factorset}
\[
f: C \x C \xra{\quad} A;\qquad (c,c') \mapsto s(c) + s(c') -
s(c+c').
\]
Our main concern is the computation of such a function $f$, via
generators of $C$ and an invariant for $A$ in the above sequence
(\ref{Thm_OrderTwo}) of abelian groups.

\begin{rem}
The Connolly--Davis computation of $NL_3(\Z) \iso
NL_4(\Z,(2)^\infty)$ involves generators $\cN_{p,g}$ indexed by
polynomials $p,g \in \Z[x]$. Either $p$ or $g$ must have zero
constant coefficient, and each generator is defined as the
nonsingular $(+1)$-quadratic linking form
\[
\cN_{p,g} := \prn{\bigoplus_2 \Z[x]/2, \PMatrix{p/2 & 1/2\\ 1/2 &
0}, \PMatrix{p/2\\ g} }
\]
of exponent two over $(\Z[x],(2)^\infty)$, see \cite[Dfn. 1.6 and p.
1057]{CD}. For our computation, we identify it with a choice of
resolution by a nonsingular split $(-1)$-quadratic formation
\[
\cN_{p,g} = \prn{\bigoplus_2 \Z[x], (\PMatrix{p & 1\\ 1 & 2g\\
2 & 0\\ 0 & 2}, \PMatrix{p & 1\\ 1 & 2g}) \bigoplus_2 \Z[x]}.
\]
\end{rem}

\begin{defn}[{\cite[p. 69]{RanickiExact}}]
Let $R$ be a ring with involution, and let $F, G$ be finitely
generated projective $R$-modules.
A \textbf{nonsingular split $\epsilon$-quadratic formation $\prn{F,(\PMatrix{\gamma\\
\mu}, \theta) G}$ over $R$} consists of the hyperbolic
$\epsilon$-quadratic form
\[
\fH_\epsilon(F) := \prn{ F \oplus F^*, \PMatrix{0 & \Id_F\\ 0 & 0} }
\]
along with the standard lagrangian $F \oplus 0$, a second lagrangian
\[
\Img( \PMatrix{\gamma\\ \mu} : G \to F \oplus F^* ),
\]
and a hessian $\theta: G \to G^*$, which is a de-symmetrization of
the pullback form:
\[
\theta - \epsilon \theta^* = \PMatrix{\gamma \\ \mu}^* \PMatrix{0 &
\Id_F\\ 0 & 0} = \gamma^* \circ \mu: G \longra G^*.
\]\qed
\end{defn}

\begin{defn}\label{Def_QM}
For any polynomial $q\in x\Z[x]$, define the nonsingular split
$(-1)$-quadratic formation $\cQ_q$ over $\Z[C_2][x]$, where $\hat{q}
:= 2(1-T)q$, by
\[
\cQ_q := \prn{\bigoplus_2 \Z[C_2][x], (\PMatrix{ 0 & \hat{q} \\ \hat{q} & 0\\
1 & (1-T)q\\(1-T) & 1}, \PMatrix{\hat{q} & 0\\ \hat{q} & q \hat{q}})
\bigoplus_2 \Z[C_2][x] }.
\]

For any polynomials $p,g \in \Z[x]$ with $pg \in x \Z[x]$, define
the nonsingular split $(-1)$-quadratic formation $\cM_{p,g}$ over
$\Z[C_2][x]$ by
\[
\cM_{p,g} := \prn{\bigoplus_2 \Z[C_2][x], (\PMatrix{p & 1\\
1 & (1-T)g\\ 2 & 0\\ 0 & 2}, \PMatrix{p & 1\\ 1 & (1-T)g})
\bigoplus_2 \Z[C_2][x]}.
 \]\qed
\end{defn}

Indeed each of these $(-1)$-quadratic formations consists of
lagrangian summands, since the associated 1-dimensional
$(-1)$-quadratic complex over $\Z[C_2][x]$ is connected \cite[Proof
2.3]{RanickiAlgebraicI} and in fact Poincar\'e: the Poincar\'e
duality map on the level of projective modules induces isomorphisms
on the homology groups. For example in $\cM_{p,g}$, the nontrivial
homological Poincar\'e duality map is
\[
\SmMatrix{p & 1\\ 1 & (1-T)g} : H^0(C) \to H_1(C),
\quad\text{where}\quad H^0(C) = H_1(C) = \Z[C_2][x]/2.
\]
Its determinant $(1-T)pg - 1$ is a unit mod $2$ in the commutative
ring $\Z[C_2][x]$, since
\[
\prn{(1-T)pg - 1}^2 = 2(1-T)(pg)^2 - 2(1-T)pg + 1 \equiv 1 \pmod{2}.
\]
Therefore the Poincar\'e duality map for $\cM_{p,g}$ is a homology
isomorphism. Also, the formation $\cQ_q$ is obtained as a pullback
of a nonsingular formation, cf. Proof \ref{Prop_QM}(1).

\begin{prop}\label{Prop_QM}
The following formulas are satisfied for cobordism classes in the
reduced module $\wt{NL}_3(\Z[C_2])$.
\begin{enumerate}
\item Boundary map: $\wt{\bdry}[P_{q,1}] = [\cQ_q]$
\item Lifts: $i^- [\cM_{p,g}] = [\cN_{p,g}]$ and $i^+ [\cM_{p,g}] =
0$
\end{enumerate}
\end{prop}

Now we state the basic relations between our generators $\cQ$ and
$\cM$, established by algebraic surgery.  Their inspiration is the
statement and proof of \cite[Lemma 4.3]{CD}, but they are proven
independently.

\begin{prop}\label{Prop_NL3_Relations}
The following formulas are satisfied for cobordism classes in the
reduced module $\wt{NL}_3(\Z[C_2])$.
\begin{enumerate}
\item Additivity: $[\cM_{p_1,g}] + [\cM_{p_2,g}] = [\cM_{p_1+p_2,g}] + [\cQ_q]$ where $q := (p_1 g) (p_2 g)$
\item Symmetry: $[\cM_{2p,g}] = [\cM_{2g,p}]$
\item Square associativity: $[\cM_{x^2 p, g}] = [\cM_{p, x^2 g}]$
\item Square root: $[\cM_{2 p^2 g, g}] = [\cM_{2 p,g}]$
\end{enumerate}
\end{prop}

Here are some useful formal consequences, which do not require the
technique of algebraic surgery.

\begin{cor}\label{Cor_NL3_Relations}
The following formulas are satisfied for cobordism classes in the
reduced module $\wt{NL}_3(\Z[C_2])$.
\begin{enumerate}
\item Exponent four: $4 \cdot [\cM_{p,g}] = 0$
\item Idempotence: $2(V_2 - 1) \cdot [\cM_{p,1}] =
0$
\item Exponent two: $2 \cdot \prn{[\cM_{x,g}] - [\cM_{1,xg}]} = 0$
\item Nilpotence: $V_2 \cdot \prn{[\cM_{x,g}] - [\cM_{1,xg}]} = 0$
\end{enumerate}
\end{cor}

Finally we conclude with a determination of the Verschiebung module
extension $\UNil_3(\Z[C_2])$, through the eyes of the
Connolly--Ranicki isomorphism (\ref{Thm_CRiso}).

\begin{thm}\label{Thm_NL3}
The extension of $\cV$-modules in Theorem \ref{Thm_OrderTwo} is
trivial.
\end{thm}

\section{Main proofs using relations}

\begin{proof}[Proof of Theorem \ref{Thm_VanishExponentTwo}]
Denote $S$ as the Sylow 2-subgroup of $F$.  Since $S$ is normal and
abelian, by the reduction isomorphism of \cite[Theorem
1.1]{Khan_Sylow2NormalAbelian} and the Connolly--Ranicki isomorphism
$r$ of Theorem \ref{Thm_CRiso}, it suffices to show that:
\[
NL_n^h(\Z[S]) = 0 \quad\text{if } n \equiv 0,1 \pmod{4}
\]
and the following induced map is an isomorphism:
\[
NL_2^h(\Z[S]) \xra{\quad} NL_2^h(\F_2).
\]

We induct on the order of $S$.  If $\abs{S}=1$, then recall from
Remark \ref{Rem_NLcomp} that
\[
NL_n(\Z[S]) = NL_n(\Z[1]) = 0 \quad\text{if } n \equiv 0,1 \pmod{4}
\]
and the following induced map is an isomorphism:
\[
NL_2(\Z[S]) = NL_2(\Z[1]) \xra{\quad} NL_2(\F_2).
\]

Otherwise suppose $\abs{S} > 1$.  Since $S$ has exponent two, there
is a decomposition
\[
S = S' \x C_2
\]
as an internal direct product of groups of exponent two.  Then the
Mayer--Vietoris sequence of \cite[Proposition
6.1]{Khan_Sylow2NormalAbelian} specializes to:
\[
\cdots NL_{n+1}(\F_2) \xra{\en\bdry\en} NL_n(\Z[S]) \xra{\quad}
\bigoplus_2 NL_n(\Z[S']) \xra{\quad} NL_n(\F_2) \xra{\en\bdry\en}
\cdots.
\]
Observe, by inductive hypothesis and Remark \ref{Rem_NLcomp}, that
$NL_n(\Z[S']) = NL_n(\F_2) = 0$ for all $n \equiv 0,1 \pmod{4}$. So
we obtain $NL_0(\Z[S]) = 0$ and an exact sequence
\[
0 \xra{\quad} NL_2(\Z[S]) \xra{\quad} \bigoplus_2 NL_2(\Z[S'])
\xra{\quad} NL_2(\F_2) \xra{\en\bdry\en} NL_1(\Z[S]) \xra{\quad} 0.
\]
But the following induced map is an isomorphism, by inductive
hypothesis:
\[
NL_2(\Z[S']) \xra{\quad} NL_2(\F_2).
\]
Therefore
\[
NL_1(\Z[S]) = 0
\]
and the following composite of induced maps is an isomorphism:
\[
NL_2(\Z[S]) \xra{\quad} NL_2(\Z[S']) \xra{\quad} NL_2(\F_2).
\]
This concludes the induction on $\abs{S}$.
\end{proof}

\begin{proof}[Proof of Theorem \ref{Thm_OrderTwo}]
The exact sequence of \cite[Proposition
6.1]{Khan_Sylow2NormalAbelian} becomes
\[
NL_{n+1}(\F_2) \xra{\en\bdry\en} NL_n(\Z[C_2]) \xra{\SmMatrix{i^-\\
i^+}} NL_n(\Z) \oplus NL_n(\Z) \xra{\SmMatrix{j^- & -j^+}}
NL_n(\F_2).
\]
Since this sequence is functorial, it must consist of $\cV$-module
morphisms. It follows by Orientable Reduction \cite[Prop.
4.1]{Khan_Sylow2NormalAbelian} that there is the commutative diagram
of Figure \ref{Fig_NLReduced} with top row exact,
\begin{figure}[!ht]
\begin{center}\tiny
$\begin{diagram} \node{0} \arrow{e} \node{NL_4(\F_2)}
\arrow{e,t}{\bdry}
\arrow[2]{s,=} \node{NL_3(\Z[C_2])} \arrow[2]{s,t}{\Id-\vare} \arrow{e,t}{\SmMatrix{i^-\\
i^+}} \arrow{sse,t}{i^- - i^+} \node{NL_3(\Z) \oplus NL_3(\Z)}
\arrow{e} \arrow[2]{s,b}{\proj_{\text{skew-diag}}} \node{0}\\
\\ \node{0} \arrow{e} \node{NL_4(\F_2)} \arrow{e,t}{\wt{\bdry}}
\node{\wt{NL}_3(\Z[C_2])} \arrow{e,t}{i^-} \node{NL_3(\Z)}
\arrow{e} \node{0}
\end{diagram}$
\end{center}
\caption{Reduction of $NL_3$}\label{Fig_NLReduced}
\end{figure}
where $\vare: C_2 \to C_2$ is the trivial map and
\[
\wt{\bdry} := (\Id - \vare) \circ \bdry.
\]
The map $\wt{\bdry}$ is a monomorphism, since $i^+\circ \bdry = 0$
and the left square commutes. The map $i^-$ is an epimorphism, since
the projection $\proj_{\text{skew-diag}}$ onto the skew-diagonal is
surjective and the right square commutes. Exactness at
$\wt{NL}_3(\Z[C_2])$ follows from its definition and exactness of
the top row at $NL_3(\Z[C_2])$.  Thus the bottom row exists and is
an exact sequence of $\cV$-modules.
\end{proof}

\begin{proof}[Proof of Proposition \ref{Prop_QM}(1)]
According to \cite[pp. 517--519]{RanickiExact}, the boundary map
\[
\bdry = \bdry_{i^-} \circ \delta: L_4(\F_2[x]) \to L_3(\Z[C_2][x])
\]
for our cartesian square is defined in general in terms of pullback
modules by
\begin{multline*}
(A'^r, \psi') \longmapsto\\ \prn{ (B^r,\Id,B'^r) , (\PMatrix{ (\Id - (\chi + \chi^*) \circ \phi, \Zero) \\
(\phi,\Id) }, (\psi - \phi \circ \chi \circ \phi, \Zero) )
(B^r,\phi', B'^r) }.
\end{multline*}
It sends a Witt class of a rank $r$ nonsingular form over $A' =
\F_2[x]$ to the Witt class of split formation over $A = \Z[C_2][x]$
obtained by pullback of the boundary formation of the lifted form
over $B = \Z[x]$ and of the hyperbolic formation over  $B' = \Z[x]$.
The form $\psi$ over $B$ lifts the input form $\psi'$ over $A'$.
Their symmetrizations are denoted
\[
\phi := \psi + \psi^*: B^r \xra{\quad} (B^r)^* \quad\text{and}\quad
\phi' := \psi' + \psi'^*: B'^r \xra{\quad} (B'^r)^*.
\]
The morphism $\chi: (B'^r)^* \to B'^r$ lifts the map
\[
\chi' := (\phi')\inv \circ \psi' \circ (\phi')\inv : (A'^r)^*
\xra{\quad} A'^r.
\]

Now we compute these morphisms in our situation.  Let $p \in x
\Z[x]$. Recall (\ref{Rem_NLcomp}) and take
\[
(A'^r,\psi') = P_{q,1} = \prn{\F_2[x]^2, \Matrix{q & 1\\
0 & 1}}.
\]
Choose a lift
\[
(B^r, \psi) = \prn{ \Z[x]^2, \Matrix{q & 1\\ 0 & 1} }.
\]
Then we obtain and select
\[
\chi' = \Matrix{1 & 0\\ 1 & q} : \F_2[x]^2 \to \F_2[x]^2
\qquad\qquad \chi = \Matrix{-1 & 0\\ 1 & -q} : \Z[x]^2 \to \Z[x]^2.
\]
Using the pullback module structure \cite[p. 507]{RanickiExact}
\[
\Z[C_2][x] \xra{\en\iso\en} (\Z[x], \Id: \F_2[x] \to \F_2[x],
\Z[x]);\qquad (m + nT) \longmapsto (m-n,m+n),
\]
the pullback formation is
\begin{align*}
\bdry[P_{q,1}] &= \prn{ (\Z[x]^2,\Id, \Z[x]^2), (\PMatrix{ (\Matrix{ 4q & 0\\ 0 & 4q }, \Zero) \\
( \Matrix{ 2q & 1\\ 1 & 2}, \Id) }, (\Matrix{ 4q^2 & 4q \\ 0 & 4q}, \Zero) ) (\Z[x]^2, \Matrix{0 & 1\\
1 & 0}, \Z[x]^2) }\\ \\
&= \prn{ (\Z[x]^2,\Id, \Z[x]^2), (\PMatrix{ (\Matrix{ 0 & 4q\\ 4q & 0 }, \Zero) \\
( \Matrix{ 1 & 2q\\ 2 & 1}, \Id) }, (\Matrix{ 4q & 0\\ 4q & 4q^2},
\Zero) ) (\Z[x]^2, \Id, \Z[x]^2) }\\ \\
&= \prn{ \Z[C_2][x]^2, (\Matrix{  0 & 2(1-T)q \\ 2(1-T)q & 0\\
1 & (1-T)q \\ (1-T) & 1}, \Matrix{ 2(1-T)q & 0\\ 2(1-T) q & 2(1-T) q^2 }) \Z[C_2][x]^2 } \\ \\
&= [\cQ_q].
\end{align*}
\end{proof}

\begin{proof}[Proof of Proposition \ref{Prop_QM}(2)]
Clearly $i^-(\cM_{p,g}) = \cN_{p,g}$. Note that the second
lagrangian $G$ of $i^+(\cM_{p,g})$ is
\[
\Img\PMatrix{p & 1\\ 1 & 0\\ 2 & 0\\ 0 & 2} = \Img\PMatrix{0 & 1\\
1 & 0\\ 2 & 0\\ -2p & 2} = \Img\PMatrix{1 & 0\\ 0 & 1\\ 0 & 2\\ 2 &
-2p}.
\]
Therefore $i^+(\cM_{p,g})$ is a graph formation over $\Z[x]$, hence
represents $0$ in $NL_3(\Z)$.
\end{proof}

\begin{proof}[Proof of Corollary \ref{Cor_NL3_Relations}(1)]
Note, by Proposition \ref{Prop_NL3_Relations}(1) and the relations
(\ref{Rem_NLcomp}) in $NL_4(\F_2)$, that
\[
4 \cdot [\cM_{p,g}] = 2 \cdot [\cQ_{pg}] + 2 \cdot [\cM_{2p,g}]
= 2 \cdot [\cQ_{pg}] + [\cQ_{2pg}] + [\cM_{4p,g}]
= [\cM_{4p,g}].
\]
There is an isomorphism
\[
(\Id, \Id, \PMatrix{p & 0\\ 0 & 0}): \cM_{0,g} \xra{\quad}
\cM_{4p,g}
\]
of split $(-1)$-quadratic formations over $\Z[C_2][x]$, see \cite[p.
69, defn.]{RanickiExact}. Therefore, as cobordism classes in
$NL_3(\Z[C_2])$, we obtain
\[
4 \cdot [\cM_{p,g}] = [\cM_{4p,g}]
= [\cM_{0,g}]\\
= 0,
\]
by \cite[Proposition 1.6.4]{RanickiExact} and since $\cM_{0,g}$ is a
graph formation.
\end{proof}

\begin{proof}[Proof of Corollary \ref{Cor_NL3_Relations}(2)]
Note, by Proposition \ref{Prop_NL3_Relations}(1) and the relations
(\ref{Rem_NLcomp}) in $NL_4(\F_2)$, that
\begin{eqnarray*}
2(V_2 - 1) \cdot [\cM_{p,1}] &=& (V_2 - 1) \cdot \prn{ [\cQ_{p}] +
[\cM_{2p,1}] }\\
&=& [\cM_{2 V_2 p,1}] - [\cM_{2p,1}]\\
&=& \sum_{k=1}^n p_k \cdot \prn{[\cM_{2 (x^k)^2,1}] -
[\cM_{2(x^k),1}]}
\end{eqnarray*}
where we write the polynomial
\[
p = p_1 x + \cdots + p_n x^n \in \Ker(\aug_0)
\]
for some $n \in \N$ and $p_1, \ldots, p_n \in \Z$. But by
Proposition \ref{Prop_NL3_Relations}(4), we have
\[
[\cM_{2 (x^k)^2,1}] - [\cM_{2(x^k),1}] = 0 \quad\text{for all } k >
0.
\]
Therefore
\[
2(V_2-1) \cdot [\cM_{p,1}] = 0.
\]
\end{proof}

\begin{proof}[Proof of Corollary \ref{Cor_NL3_Relations}(3)]
Note by Proposition \ref{Prop_NL3_Relations}(1,2) and the relations
(\ref{Rem_NLcomp}) in $NL_4(\F_2)$ that
\begin{eqnarray*}
2\cdot ([\cM_{x,g}] - [\cM_{1,xg}]) &=& ([\cQ_{xg}] - [\cQ_{xg}]) +
([\cM_{2x,g}] - [\cM_{2,xg}])\\
&=& [\cM_{2g,x}] - [\cM_{2xg,1}].
\end{eqnarray*}
By Proposition \ref{Prop_NL3_Relations}(1) using the fact that
$[\cQ_{q}] = 0$ if $q$ is a multiple of $2$, and since $g$ has
$\Z$-coefficients, we may assume that $g = x^k$ for some $k \in \N$
in order to show that the right-hand term vanishes. If $k=2i$ is
even, then by Proposition \ref{Prop_NL3_Relations}(3,2), note
\[
[\cM_{2g,x}] = [\cM_{2(x^i)^2,x}] = [\cM_{2,(x^i)^2 x}] =
[\cM_{2,xg}] = [\cM_{2xg,1}].
\]
Otherwise suppose $k=2i+1$ is odd. Then by Proposition
\ref{Prop_NL3_Relations}(4) twice and by induction on $k$, note
\[
[\cM_{2g,x}] = [\cM_{2(x^i)^2 x, x}]
= [\cM_{2(x^i),x}]
= [\cM_{2 (x^{i+1}),1}]
= [\cM_{2 (x^{i+1})^2,1}]
= [\cM_{2 xg,1}].
\]
Therefore for all $g \in \Z[x]$ we obtain
\[
2\cdot ([\cM_{x,g}] - [\cM_{1,xg}]) = 0.
\]
\end{proof}

\begin{proof}[Proof of Corollary \ref{Cor_NL3_Relations}(4)]
Note by Proposition \ref{Prop_NL3_Relations}(3) that
\[
V_2 \cdot \prn{[\cM_{x,g}] - [\cM_{1,xg}]} = [\cM_{x^2,V_2 g}] -
[\cM_{1, x^2 V_2 g}] = 0.
\]
\end{proof}

\begin{proof}[Proof of Theorem \ref{Thm_NL3}]
Consider the $\cV$-module morphism
\[
s: NL_3(\Z) \xra{\quad} \wt{NL}_3(\Z[C_2])
\]
given additively by
\begin{gather*}
V_n \cdot [\cN_{x,1}] \longmapsto V_n \cdot [\cM_{x,1}]\\
V_n \cdot ([\cN_{x,x^p}] - [\cN_{1,x x^p}]) \longmapsto V_n \cdot
([\cM_{x,x^p}] - [\cM_{1,x x^p}]).
\end{gather*}
The section $s$ is a well-defined $\cV$-module morphism by Corollary
\ref{Cor_NL3_Relations}.
\end{proof}

\section{Some algebraic surgery machines}

The remaining proofs of all parts of Proposition
\ref{Prop_NL3_Relations} are technical---algebraic surgery and
gluing are required. The first machine has input certain quadratic
formations and has output quadratic forms.

\begin{lem}\label{Lem_BoundaryInstant}
Suppose $(C,\psi)$ is a 1-dimensional $(-1)$-quadratic Poincar\'e
complex over $\Z[C_2][x]$ satisfying the following hypotheses.
\begin{enumerate}
\item[(a)] The 1-dimensional chain complex $C$ over $\Z[C_2][x]$ has modules $C_1 =
C_0$ and differential $d_C = 2 \cdot \Id$.
\item[(b)] There is a null-cobordism
\[
(f: i^-(C) \to D, (\delta \psi, i^-(\psi)) \in W_\%(f,-1)_2)
\]
over $\Z[x]$ such that $f_0 = \Id: C_0 \to D_0$ and $\delta \psi_2 =
0: D^0 \to D_0$.
\item[(c)] The quadratic Poincar\'e complex
$i^+(C,\psi)$ over $\Z[x]$ corresponds to a graph formation.
\end{enumerate}

Then we obtain the following conclusions.
\begin{enumerate}
\item
There exists a 2-dimensional $(-1)$-quadratic Poincar\'e complex
$(F,\Psi)$ over $\F_2[x]$ such that
\[
\bdry[\ol{S}(F,\Psi)] = [\ol{S}(C,\psi)].
\]
Here,
\[
\bdry: L_4(\F_2[x]) \xra{\quad} L_3(\Z[C_2][x])
\]
is the boundary map of the Mayer--Vietoris sequence of Rim's
cartesian square, and $\ol{S}$ is the skew-suspension
isomorphism.

\item The instant surgery obstruction $\Omega(F,\Psi)$ is Witt
equivalent to the nonsingular $(+1)$-quadratic form $j^-(D^1,
\delta\psi_0)$ over $\F_2[x]$.
\end{enumerate}
\end{lem}

The next machine constructs inputs for the above one given a
lagrangian of a certain linking form. It is obtained as a
specialization of \cite[Proof 3.4.5(ii)]{RanickiExact}\footnote{See
errata for the formulas at
\texttt{http://www.maths.ed.ac.uk/\~{}aar/books/exacterr.pdf}.}.

\begin{lem}\label{Lem_NullCobordism}
Suppose $(C,\psi)$ is a 1-dimensional $(-1)$-quadratic Poincar\'e
complex over $\Z[C_2][x]$ satisfying the following hypotheses.
\begin{enumerate}
\item[(a)] The 1-dimensional chain complex $C$ over $\Z[C_2][x]$ has modules $C_1 =
C_0$ and differential $d_C = 2 \cdot \Id$.

\item[(b)] There exists a lagrangian $L$ of the nonsingular $(+1)$-quadratic linking form $(N,b,q)$
 over $(\Z[x], (2)^\infty)$ associated to $i^-(C,\psi)$.

\item[(c)] The evaluation $i^+(C,\psi)$ corresponds to a graph formation over $\Z[x]$.
\end{enumerate}

Choose a finitely generated projective module $P$ over $\Z[C_2][x]$
and morphisms
\begin{enumerate}
\item[(i)] $\pi: P \to C^1$ monic with image
\[
i^-(\pi)(P) = e\inv(L),
\]
where the quotient map $e$ is
\[ e : i^-(C^1) \xra{\quad} N :=
Cok\prn{i^-(d_C^* : C^0 \to C^1)}
\]

\item[(ii)] $\chi: P \to P^*$ satisfying the de-symmetrization
identity
\[
(\pi\inv \circ d_C^*)^* \circ (\chi + \chi^*) = (\wt{\psi}_0
 - \psi_0^*) \circ \pi : P \to C_0.
\]
\end{enumerate}

Then we obtain the following conclusions.
\begin{enumerate}
\item Define a quadratic cycle $\wh{\psi} \in W_\%(C,-1)_1$ by
\begin{gather*}
\wh{\psi}_0 := \psi_0: C^0 \xra{\quad} C_1 \qquad\qquad
\wt{\wh{\psi}}_0 := \wt{\psi}_0 : C^1 \xra{\quad} C_0
\\
\wh{\psi}_1 := (\pi\inv \circ d_C^*)^* \circ \chi \circ (\pi\inv
\circ d_C^*) - \wt{\psi}_0 \circ d_C^* : C^0 \xra{\quad} C_0.
\end{gather*}
Then the quadratic cycle $\wh{\psi}$ is homologous to $\psi$ in
$W_\%(C,-1)_1$ over $\Z[C_2][x]$.

\item Define a chain complex $D = \set{D_1 \xra{d_D} D_0}$ with modules
\[
D_1 := i^-(P^*) \quad\text{and}\quad D_0 := i^-(C_0)
\]
and with differential
\[
d_D := i^-(\pi\inv \circ d_C^*)^*.
\]
Define a chain map $f: i^-(C) \to D$ by
\[
f_0 := \Id: i^-(C_0) \xra{\quad} D_0 \qquad\qquad f_1 := i^-(\pi^*)
: i^-(C_1) \xra{\quad} D_1.
\]
Define a quadratic chain $\delta\psi \in W_\%(D,-1)_2$ by
\begin{gather*}
\delta\psi_0 := - i^-(\chi)^* : D^1 \xra{\quad} D_1\qquad\qquad
\delta\psi_1 := -i^-(\chi) \circ d_D^*: D^0 \xra{\quad} D_1
\\
\wt{\delta\psi}_1 := \wt{\psi}_0 \circ i^-(\pi): D^1 \xra{\quad}
D_0\qquad\qquad \delta\psi_2 := 0: D^0 \xra{\quad} D_0.
\end{gather*}
Then
\[
\prn{f: i^-(C) \to D, (\delta\psi, i^-(\wh{\psi})) \in
W_\%(f,-1)_2}\] is a null-cobordism over $\Z[x]$.

\item The evaluation $i^+(C,\wh{\psi})$ corresponds to a graph
formation over $\Z[x]$.
\end{enumerate}
\end{lem}

Composition of the lemmas yields immediately the following result.

\begin{prop}\label{Prop_NullInstant}
Suppose $(C,\psi)$ satisfies Hypotheses (a,b,c) of Lemma
\ref{Lem_NullCobordism}, and choose $P, \pi, \chi$ accordingly. Then
there exists a 2-dimensional $(-1)$-quadratic complex $(F,\Psi)$
over $\F_2[x]$ such that
\[
\bdry[(F,\Psi)] = [(C,\wh{\psi})] = [(C,\psi)]
\]
as cobordism classes in
\[
L_1(\Z[C_2][x],-1) \xra{\en\ol{S} \iso\en} L_3(\Z[C_2][x]),
\]
and that its instant surgery obstruction is
\[
[\Omega(F,\Psi)] = [j^-(D^1,\delta\psi_0)] = [k(P,-\chi^*)]
\]
as Witt classes of the nonsingular $(+1)$-quadratic forms over
$\F_2[x]$. \qed
\end{prop}

Now we show why these machines work.

\begin{proof}[Proof of Lemma \ref{Lem_BoundaryInstant}]
We shall put together the information in the hypotheses using a
technique called ``algebraic gluing'' \cite[\S 1.7]{RanickiExact}.
The resultant object $(F,\Psi)$ is a union \cite[pp.
77--78]{RanickiExact} over $\F_2[x]$. The more efficient ``direct
union'' \cite[pp. 79--80]{RanickiExact} does not apply here since
the null-cobordisms $(D,\delta\psi)$ and $(E,0)$ are non-split in
general.

First, define a chain complex $E = \set{E_1 \to 0}$ over $\Z[x]$
with module
\[
E_1 := i^+(C_1),
\]
and a chain map $g: i^+(C) \to E$ by
\[
g_1 := \Id: i^+(C_1) \to E_1.
\]
Then the quadratic pair
\[
\prn{g: i^+(C) \to E, (0,i^+(\psi) \in W_\%(g,-1)_2)}
\]
is the data for an algebraic surgery. Consider the 2-dimensional
mapping cone
\[
\fC(g) = \prn{ i^+(C_1) \xra{\SmMatrix{-\Id\\ 2 \cdot \Id}} E_1
\oplus i^+(C_0) \xra{\qquad} 0 }.
\]
Note that
\[
H^2(E) = H_0(\fC(g)) = 0 \quad\text{and}\quad H^0(E) = H_2(\fC(g)) =
0.
\]
Observe that
\[
H^1(E) = E^1 \quad\text{and}\quad \proj_*: H_1(\fC(g))
\xra{\en\iso\en} i^+(C_0).
\]
Then the homological Poincar\'e duality map $H^1(E) \to H_1(\fC(g))$
is given by
\[
i^+(\wt{\psi}_0 - \psi_0^*) : E^1 \xra{\quad} i^+(C_0).
\]
By hypothesis, $i^+(C,\psi)$ represents a graph formation
\[
\prn{F, (\PMatrix{\gamma\\ \mu}, \theta) G},
\]
which means that $\gamma: G \to F$ an isomorphism. According to
\cite[Proof 2.5]{RanickiAlgebraicI}, the representation is given by
\begin{gather*}
F = i^+(C_1) \quad\text{and}\quad G = i^+(C^0)
\\
\gamma = i^+(\wt{\psi}_0^* - \psi_0) \quad\text{and}\quad \mu =
i^+(d_C^*) \quad\text{and}\quad \theta = -i^+(\psi + d_C \circ
\psi_0).
\end{gather*}
Thus the map $H^1(E) \to H_1(\fC(g))$ is given by the isomorphism
$\gamma^*$.
 Since the Poincar\'e duality map $E^{2-*} \to
\fC(g)$ of projective module chain complexes induces isomorphisms in
homology, it must be a chain homotopy equivalence. Thus the
following 2-dimensional $(-1)$-quadratic pair is Poincar\'e:
\[
\prn{g: i^+(C) \to E, (0, i^+(\psi))}.
\]

Next, define a 2-dimensional $(-1)$-quadratic Poincar\'e complex
$(F,\Psi)$ over $\F_2[x]$ as the union (see \cite[pp.
77--78]{RanickiExact})
\[
(F, \Psi) := j^-\prn{f: i^-(C) \to D, (-\delta\psi,-i^-(\psi))}
\bigcup_{k(C,\psi)} j^+\prn{g: i^+(C) \to E, (0, i^+(\psi))},
\]
where $k$ is composite morphism of rings with involution:
\[
k := j^- \circ i^- = j^+ \circ i^+ : \Z[C_2] \xra{\quad} \F_2.
\]
By construction,
\[
\bdry[(F,\Psi)] = [(C,\psi)],
\]
where the boundary map
\[
\bdry: L_4(\F_2[x]) \xra{\quad} L_3(\Z[C_2][x])
\]
is defined in \cite[Props. 6.3.1, 6.1.3]{RanickiExact} for our
cartesian square. For simplicity, we suppress the morphisms $i^\pm,
j^\pm, k$ in the remainder of the proof.

The 2-dimensional chain complex $F$ over $\F_2[x]$ has modules
\[
F_2 = C_1 \qquad\qquad F_1 = D_1 \oplus C_0 \oplus E_1 \qquad\qquad
F_0 = D_0
\]
and differentials
\[
d_F^2 = \PMatrix{-f_1\\ d_C\\ -g_1}: F_2 \to F_1 \qquad\qquad d_F^1
= \PMatrix{d_D & f_0 & 0} : F_1 \to F_0.
\]

The quadratic cycle $\Psi \in W_\%(F,-1)_2$ has components
\begin{gather*}
\Psi_0^2 = \PMatrix{-\psi_0 \circ f_0^*}: F^0 \to F_2 \qquad
\Psi_0^1 = \PMatrix{-\delta\psi_0 & 0 & 0\\ \wt{\psi}_0 \circ f_1^*
& \psi_1^* & 0\\ 0 & g_1 \circ \psi_0 & 0}: F^1 \to F_1
\\
\Psi_0^0 = \PMatrix{0}: F^2 \to F_0 \qquad \Psi_1^1 = \PMatrix{-\delta\psi_1\\
\psi_1 \circ f_0^*\\ 0}: F^0 \to F_1
\\
\Psi_1^0 = \PMatrix{-\widetilde{\delta\psi}_1 & 0 & 0}: F^1 \to F_0
\qquad \Psi_2^0 = \PMatrix{0}: F^0 \to F_0.
\end{gather*}

The differential
\[
(d_F^1)^*: F^0 \xra{\quad} F^1
\]
is a split monomorphism, since $f_0 = \Id: C_0 \to D_0$. Hence the
instant surgery obstruction \cite[Prop. 4.3]{RanickiAlgebraicI} is
represented by
\[
\Omega(F,\Psi) = \prn{D^1 \oplus E^1 \oplus C_1, \PMatrix{\delta\psi_0 & 0 & -f_1\\
0 & 0 & -\Id\\ 0 & 0 & 0}}. \] This is Witt equivalent to the
(necessarily) nonsingular (+1)-quadratic form $(D^1,\delta\psi_0)$
over $\F_2[x]$.
\end{proof}

\begin{proof}[Proof of Lemma \ref{Lem_NullCobordism}]
Indeed $\wh{\psi} \in W_\%(C,-1)_1$ is a quadratic cycle, since
\begin{align*}
\wh{\psi}_1 + \wh{\psi}_1^* &= (\pi\inv \circ d_C^*)^* \circ (\chi +
\chi^*) \circ (\pi\inv \circ d_C^*) - (\wt{\psi}_0 \circ d_C^* + d_C
\circ \wt{\psi}_0^*)\\
&= (\wt{\psi}_0 - \psi_0^*) \circ d_C^* - \wt{\psi}_0 \circ d_C^* -
d_C \circ \wt{\psi}_0^*\\
&= -(d_C \circ \psi_0 + \wt{\psi}_0 \circ d_C^*)^*\\
&= (\psi_1 + \psi_1^*)^* \\
&= \psi_1 + \psi_1^*.
\end{align*}

A similar check shows that $f: i^-(C) \to D$ is a chain map and that
\[
\prn{f: i^-(C) \to D, (\delta\psi, i^-(\wh{\psi}))}
\]
is a 2-dimensional $(-1)$-quadratic pair over $\Z[x]$.  It is
Poincar\'e (see \cite[p. 259]{RanickiExact}), since it is the data
for an algebraic surgery to a contractible complex, killing the lift
$i^-(P)$ of the lagrangian $L$.

The quadratic cycles $\wh{\psi}$ and $\psi$ are
homologous\footnote{In general, $\wh{\psi}$ and $\psi$ are
\emph{quadratic} homotopy equivalent \cite[p. 71 Defn., Prop.
3.4.5(ii)]{RanickiExact}.}: the differences $\wh{\psi}_0 - \psi_0$
and $\wt{\wh{\psi}}_0 - \wt{\psi}_0$ are zero, and the difference
$\wh{\psi}_1 - \psi_1$ is $(-1)$-symmetric (see above calculation).
Therefore, the latter difference is $(-1)$-even since
$\wh{H}^0(\Z[x],-1) = 0$. Finally, $i^+(C,\wh{\psi})$ corresponds to
the same graph formation as $i^+(C,\psi)$, except that their
hessians have difference $\wh{\theta} - \theta = \psi_1 -
\wh{\psi}_1$.
\end{proof}

\section{Remaining proofs of relations}

Using our machine (\ref{Prop_NullInstant}), we grind out the primary
relations (\ref{Prop_NL3_Relations}) in $\wt{NL}_3(\Z[C_2])$ as a
$\cV$-module.

\begin{proof}[Proof of Proposition \ref{Prop_NL3_Relations}(1)]
Let $(C,\psi)$ be a 1-dimensional $(-1)$-quadratic Poincar\'e
complex associated to the following nonsingular split
$(-1)$-quadratic formation over $\Z[C_2][x]$:
\[
\cM_{p_1,g} \oplus \cM_{p_2,g} \oplus \cM_{p_1+p_2,g}.
\]
In particular, it has modules $C_1 = C_0$ of rank $6$ and
differential $d_C = 2\cdot \Id$. Consider the exponent two linking
form $(N,b,q)$ over $(\Z[x],(2)^\infty)$ associated to the
evaluation $i^-(C,\psi)$, defined as
\[
(N,b,q) = \cN_{p_1,g} \oplus \cN_{p_2,g} \oplus -\cN_{p_1+p_2,g}.
\]
Define a lift $\pi: P \to C^1$ of a lagrangian $L$ of $(N,b,q)$ and
a morphism $\chi: P \to P^*$ by
\[
\pi := \Matrix{0 & 1 & 0 & 2 & 0 & 0\\ 1 & 0 & 1 & 0 & 2 & 0\\ 0 & 1
& 0 & 0 & 0 & 2\\ 1 & 0 & 0 & 0 & 0 & 0\\ 0 & 1 & 0 & 0 & 0 & 0\\
0 & 0 & 1 & 0 & 0 & 0},\quad \chi := \Matrix{g & 1 & g & 1 & 2g & 1\\ 0 & 0 & 0 & p_1 & 1 & p_2\\
0 & 0 & 0 & 1 & 2g & 0\\ 0 & 0 & 0 & p_1 & 2 & 0\\ 0 & 0 & 0 & 0 &
2g & 0\\ 0 & 0 & 0 & 0 & 0 & p_2}.
\]
It is straightforward to verify that $i^-(\pi)(P)$ is the inverse
image of a lagrangian $L$ and that $\chi$ satisfies the
de-symmetrization identity. Therefore, by Proposition
\ref{Prop_NullInstant}, we obtain a 2-dimensional $(-1)$-quadratic
Poincar\'e complex $(F,\Psi)$ over $\F_2[x]$ such that
\[
\bdry[(F,\Psi)] = [(C,\wh{\psi})] = [(C,\psi)]
\qquad\text{and}\qquad [\Omega(F,\Psi)] = [k(P, -\chi^*)].
\]

In classical notation, we have that $(F,\Psi)$ is represented by the
nonsingular $(+1)$-quadratic form
\[
(M,\lambda,\mu) := \prn{ \bigoplus_6 \F_2[x], \Matrix{ 0 & 1 & g & 1 & 0 & 1 \\
1 & 0 & 0 & p_1 & 1 & p_2 \\ g & 0 & 0 & 1 & 0 & 0 \\ 1 & p_1 & 1 & 0 & 0 & 0 \\ 0 & 1 & 0 & 0 & 0 & 0 \\
1 & p_2 & 0 & 0 & 0 & 0 }, \Matrix{ g \\ 0 \\ 0 \\ p_1 \\ 0 \\ p_2 }
}.
\]
Its pullback along the choice (see \cite[Proof 5.3]{Wall}) of
automorphism
\[
\alpha := \Matrix{ 1 & 0 & 0 & 0 & 1 & 0 \\ 0 & 1 & g & 1 & 0 & 1 \\
0 & 0 & 1 & 0 & 1 & 0 \\ 0 & 0 & 0 & 1 & g & 0 \\ 0 & 0 & 0 & p_1 & 1 + p_1 g & p_2 \\
0 & 0 & 0 & 0 & 0 & 1 } : M \xra{\quad} M
\]
is the symplectic form
\[
\alpha^*(M,\lambda,\mu) = \prn{ \bigoplus_6 \F_2[x], \Matrix{ 0 & 1
& 0 & 0 & 0 & 0 \\ 1 & 0 & 0 & 0 & 0 & 0 \\ 0 & 0 & 0 & 1 & 0 & 0 \\ 0 & 0 & 1 & 0 & 0 & 0 \\ 0 & 0 & 0 & 0 & 0 & 1\\
0 & 0 & 0 & 0 & 1 & 0 }, \Matrix{g \\ 0 \\ 0 \\ p_1 \\ p_1 g^2 \\
p_2 } }
\]
which has Arf invariant $[q]$, where $q := (p_1 g) (p_2 g)$. So, by
Remark \ref{Rem_NLcomp}, as cobordism classes in $L_4(\F_2[x])$, we
must have $[(F,\Psi)] = [P_{q,1}]$. Therefore, by Proposition
\ref{Prop_QM}(1), we obtain
\[
[\cM_{p_1,g}] + [\cM_{p_2,g}] - [\cM_{p_1+p_2,g}] = [(C,\psi)] =
\bdry[(F,\Psi)] = \bdry[P_{q,1}]  = \wt{\bdry}[P_{q,1}]
= [\cQ_{q}].
\]
\end{proof}

\begin{proof}[Proof of Proposition \ref{Prop_NL3_Relations}(2)]
Let $(C,\psi)$ be a 1-dimensional $(-1)$-quadratic Poincar\'e
complex associated to the following nonsingular split
$(-1)$-quadratic formation over $\Z[C_2][x]$:
\[
\cM_{2 p, g} \oplus -\cM_{2g, p}.
\]
In particular, it has modules $C_1 = C_0$ of rank $4$ and
differential $d_C = 2\cdot \Id$. Consider the exponent two linking
form $(N,b,q)$ over $(\Z[x],(2)^\infty)$ associated to evaluation
$i^-(C,\psi)$, defined as
\[
(N,b,q) = \cN_{2 p, g} \oplus -\cN_{2g, p}.
\]
Define a lift $\pi: P \to C^1$ of a lagrangian $L$ of $(N,b,q)$ and
a morphism $\chi: P \to P^*$ by
\[
\pi := \Matrix{1 & 0 & 2 & 0\\ 0 & 1 & 0 & 2\\ 0 & 1 & 0 & 0\\ 1 & 0
& 0 & 0},\quad \chi := \Matrix{0 & 0 & 2p & 1\\ 0 & 0 & 1 & 2g\\ 0 & 0 & 2p & 2\\
0 & 0 & 0 & 2g}.
\]
The remainder follows by an argument similar to Proof
\ref{Prop_NL3_Relations}(1).
\end{proof}

\begin{proof}[Proof of Proposition \ref{Prop_NL3_Relations}(3)]
Let $(C,\psi)$ be a 1-dimensional $(-1)$-quadratic Poincar\'e
complex associated to the following nonsingular split
$(-1)$-quadratic formation over $\Z[C_2][x]$:
\[
\cM_{x^2 p, g} \oplus -\cM_{p, x^2 g}.
\]
In particular, it has modules $C_1 = C_0$ of rank $4$ and
differential $d_C = 2\cdot \Id$. Consider the exponent two linking
form $(N,b,q)$ over $(\Z[x],(2)^\infty)$ associated to the
evaluation $i^-(C,\psi)$, defined as
\[
(N,b,q) = \cN_{x^2 p,g} \oplus -\cN_{p,x^2 g}.
\]
Define a lift $\pi: P \to C^1$ of a lagrangian $L$ of $(N,b,q)$ and
a morphism $\chi: P \to P^*$ by
\[
\pi := \Matrix{1 & 0 & 0 & 0\\ 0 & x & 0 & 2\\ x & 0 & 2 & 0\\ 0 & 1
& 0 & 0}, \quad \chi := \Matrix{0 & 0 & -xp & 1\\ 0 & 0 & -1 & 2xg\\ 0 & 0 & -p &
0\\ 0 & 0 & 0 & 2g}.
\]
The remainder follows by an argument similar to Proof
\ref{Prop_NL3_Relations}(1).
\end{proof}

\begin{proof}[Proof of Proposition \ref{Prop_NL3_Relations}(4)]
Let $(C,\psi)$ be a 1-dimensional $(-1)$-quadratic Poincar\'e
complex associated to the following nonsingular split
$(-1)$-quadratic formation over $\Z[C_2][x]$:
\[
\cM_{2 p^2 g, g} \oplus -\cM_{2p, g}.
\]
In particular, it has modules $C_1 = C_0$ of rank $4$ and
differential $d_C = 2\cdot \Id$. Consider the exponent two linking
form $(N,b,q)$ over $(\Z[x],(2)^\infty)$ associated to the
evaluation $i^-(C,\psi)$, defined as
\[
(N,b,q) = \cN_{2 p^2 g,g} \oplus -\cN_{2p, g}.
\]
Define a lift $\pi: P \to C^1$ of a lagrangian $L$ of $(N,b,q)$ and
a morphism $\chi: P \to P^*$ by
\[
\pi := \Matrix{1 & 1 & 0 & 0\\ 0 & p & 2 & 0\\ 0 & 1 & 0 & 0\\ p & 0
& 0 & 2},\quad \chi := \Matrix{0 & p^2 g & 1 & -2 pg\\ 0 & p^2 g & 1+2pg & -1\\ 0 &
0 & 2g & 0\\ 0 & 0 & 0 & -2g}.
\]
The remainder follows by an argument similar to Proof
\ref{Prop_NL3_Relations}(1).
\end{proof}

This concludes the calculation of $\UNil_*(\Z[C_2])$ as a
Verschiebung module.

\subsection*{Acknowledgements}

The author would like to thank Jim Davis for suggesting an extension
of the methodology of Connolly--Davis \cite{CD}. The results formed
a portion of the author's thesis \cite{Khan_Dissertation} under his
supervision at Indiana University.


\bibliographystyle{alpha}
\bibliography{PostDissertation}

\end{document}